\documentclass[12pt]{article}

\usepackage{latexsym, amssymb, amsmath, amscd, amsfonts, epsfig, graphicx, colordvi,verbatim,ifpdf}
\usepackage{amsfonts, amsmath, amssymb}
\usepackage{amssymb,amsfonts,amsmath,latexsym,epsfig,cite, psfrag,eepic,color}
\usepackage{amscd,graphics}
\usepackage{latexsym, amssymb,  amsmath,amscd, amsfonts, epsfig, graphicx, colordvi,amsthm}

\usepackage{graphicx}
\usepackage{color}
\usepackage{ifpdf}
\usepackage{fancybox}

\newtheorem{thm}{Theorem}[section]

\newtheorem{lem}[thm]{Lemma}

\def\pf{\noindent{\it Proof.} }
\def\qed{\nopagebreak\hfill{\rule{4pt}{7pt}}
\medbreak}

\numberwithin{equation}{section}

\def\qed{\nopagebreak\hfill{\rule{4pt}{7pt}}
\medbreak}

\newlength{\boxedparwidth}
  {\begin{center} \begin{tabular}{|@{\hspace{.315in}}c@{\hspace{.15in}}|}
                  \hline \\ \begin{minipage}[t]{\boxedparwidth}
                  \setlength{\parindent}{.25in}}%
  {\end{minipage} \\ \\ \hline \end{tabular} \end{center}}

\begin{document}
\begin{center}

{ \large\bf The Bailey transform and Hecke-Rogers identities for\\[5pt]
 the universal mock theta functions }
\end{center}

\vskip 5mm

\begin{center}
 {  Kathy Q. Ji}$^{1}$ \quad
and  \quad { Aviva X.H. Zhao}$^{2}$ \vskip 2mm

     Center for Combinatorics, LPMC-TJKLC\\[2pt]
   Nankai University, Tianjin 300071, P.R. China

   \vskip 2mm

 Email:  $^1$ji@nankai.edu.cn, \quad $^2$zhaoxiaohua@mail.nankai.edu.cn
\end{center}

\vskip 6mm \noindent {\bf Abstract.} Recently, Garvan obtained two-variable Hecke-Rogers identities for three universal mock theta functions $g_2(z;q),\,g_3(z;q),\,K(z;q)$ by using   basic hypergeometric functions,  and he proposed a problem of finding  direct proofs of these  identities by using Bailey pair technology. In this paper, we give proofs of Garvan's identities  by applying Bailey's transform with the conjugate Bailey pair of Warnaar   and three Bailey pairs deduced from two special cases of $_6\psi_6$ given by Slater. In particular, we obtain a compact form of two-variable Hecke-Rogers identity related to $g_3(z;q)$, which imply the corresponding identity given by Garvan. We also extend these two-variable Hecke-Rogers  identities into infinite families.



\section{Introduction}

This paper is concerned with two-variable Hecke-Rogers identities for the universal mock theta functions due to Garvan \cite{Garvan-2014}. Recall that universal mock theta functions are defined by  Gordon and McIntosh \cite{Gordon-2012} and Hickerson \cite{Hickerson-1988}  as
\begin{align*}
g_2(z;q)&=  \sum_{n=0}^{\infty}\frac{(-q;q)_nq^{\frac{1}{2}n(n+1)}}
  {(z;q)_{n+1}(z^{-1}q;q)_{n+1}},\\[5pt]
g_3(z;q)&=\sum_{n=0}^{\infty}
  \frac{q^{n(n+1)}}{(z;q)_{n+1}(z^{-1}q;q)_{n+1}},\\[5pt]
K(z;q) &=\sum_{n=0}^{\infty}\frac{(-1)^n(q;q^2)_nq^{n^2}}
  {(zq^2;q^2)_n(z^{-1}q^2;q^2)_n}.
\end{align*}
Here and throughout the paper,   we  adopt the standard notation on $q$-series \cite{Andrews-1976, Gasper-2004}.
\[
(a;q)_\infty =\prod_{n=0}^\infty (1-aq^n),
\]
\[ (a;q)_n=\frac{(a;q)_\infty}{(aq^n;q)_\infty},\\[3pt]
\]
\[(a_1,a_2,\ldots, a_k;q)_n=(a_1;q)_n(a_2;q)_n\cdots(a_k;q)_n,\]
and
\[(a_1,a_2,\ldots, a_k;q)_\infty=(a_1;q)_\infty(a_2;q)_\infty\cdots(a_k;q)_\infty.
\]

  Gordon and McIntosh \cite{Gordon-2012} and Hickerson \cite{Hickerson-1988}   noticed that all of Ramanujan's mock theta functions can be written in terms of  $g_2(z;q)$ and $g_3(z;q)$, that is why these two functions are called universal mock theta functions. The function $K(z;q)$ first appears in Ramanujan's lost notebook \cite{Andrews-1979} and is related by modular transformation to the universal mock theta function $g_2(z;q)$ in \cite{McIntosh-2012}.

In  \cite{Garvan-2014}, Garvan first
 obtained four Hecke-Rogers identities for these three universal mock theta functions, which  can be viewed as  two-variable generalizations of Hecke-Rogers identities.  These identities   can also be used to establish  Hecke-Rogers-type double sum representation for each of Ramanujan's mock theta functions.

More specifically,  Garvan \cite{Garvan-2014} considered the Dyson rank function $R(z;q)$ defined by
  \[R(z;q)=\sum_{n=0}^{\infty}
 \frac{q^{n^2}}{(zq;q)_n(z^{-1}q;q)_n},
  \]
  which has the following relation with $g_3(z;q)$.
   \[g_3(z;q)=-\frac{1}{z}+\frac{1 }{z(1-z)}R(z;q).\]
He obtained  the following Hecke-Rogers  identity  for $R(z;q)$.
\begin{align}\label{Garvan1}\nonumber
  &(zq, z^{-1}q,q;q)_{\infty}R(z;q)\\[3pt]\nonumber
  &=(zq, z^{-1}q,q;q)_{\infty}\sum_{n=0}^{\infty}
 \frac{q^{n^2}}{(zq;q)_n(z^{-1}q;q)_n}\\[3pt]
    &=\frac{1}{2}\sum_{n=0}^{\infty}\left(\sum_{j=0}^{[n/2]}(-1)^{n+j}(z^{n-3j}
  +z^{3j-n})q^{\frac{1}{2}(n^2-3j^2)+\frac{1}{2}(n-j)}\right.\nonumber \\
  &\quad\hskip 2cm\left.+\sum_{j=1}^{[n/2]}(-1)^{n+j}(z^{n-3j+1}
  +z^{3j-n-1})q^{\frac{1}{2}(n^2-3j^2)+\frac{1}{2}(n+j)}\right).
  \end{align}
  Putting $z=1$ in  \eqref{Garvan1},    we   see that  the Hecke-Rogers identity \cite{Hecke-1959,Rogers-1894} follows by noting  $R(1;q)=1/(q;q)_\infty$.
  \begin{equation}
  (q;q)_\infty^2=\sum_{n=0}^\infty \sum_{j=-[n/2]}^{[n/2]}(-1)^{n+j}q^{\frac{1}{2}(n^2-3j^2)+\frac{1}{2}(n+j)}.
  \end{equation}
 Garvan \cite{Garvan-2014} also considered the overpartition rank function $H(z;q)$ defined by
 \begin{equation}H(z;q)=  \sum_{n=0}^{\infty}\frac{(-1;q)_nq^{\frac{1}{2}n(n+1)}}
  {(zq;q)_n(z^{-1}q;q)_n},
 \end{equation}
 which has the following relation with $g_2(z;q)$ \cite{McIntosh-2012}.
 \begin{equation}\label{r-h-g2}
 (1+z)H(z;q)=(1-z)+2z(1-z)g_2(z;q).
 \end{equation}
 The following Hecke-Rogers identity for $H(z;q)$ is given by Garvan \cite{Garvan-2014}.
  \begin{align}\nonumber
  &(1+z)(zq, z^{-1}q,q;q)_{\infty}
 H(z;q)\\[3pt]\nonumber
 &\hskip 0.5cm=(1+z)(zq, z^{-1}q,q;q)_{\infty}
  \sum_{n=0}^{\infty}\frac{(-1;q)_nq^{\frac{1}{2}n(n+1)}}
  {(zq;q)_n(z^{-1}q;q)_n}\\[3pt]
  &\hskip 0.5cm =\sum_{n=0}^{\infty}\sum_{\mid m\mid\leqslant[n/2]}
  (-1)^{n+m}(z^{n-2|m|+1}+z^{2|m|-n})
  q^{\frac{1}{2}(n^2-2m^2)+\frac{1}{2}n} \label{Garvan2}\\[3pt]
  &\hskip 0.5cm =\sum_{n=0}^{\infty}\sum_{\mid m\mid\leqslant[n/3]}
  (-1)^{n}(z^{n-4|m|+1}+z^{4|m|-n})
  q^{\frac{1}{2}(n^2-8m^2)+\frac{1}{2}n}. \label{Garvan3}
  \end{align}
  Letting $z=1$ in \eqref{Garvan2} and \eqref{Garvan3}, and using the fact that $H(1;q)=(q^2;q^2)_\infty/(q;q)^2_\infty$,   we  get  the following two Hecke-Rogers type identities firstly given by  Kac and Peterson \cite{Andrews-1984, Kac-Peterson-1980}.
  \begin{align}
  (q;q)_\infty(q^2;q^2)_\infty& =\sum_{n=0}^{\infty}\sum_{\mid m\mid\leqslant[n/2]}
  (-1)^{n+m}q^{\frac{1}{2}(n^2-2m^2)+\frac{1}{2}n}, \\[3pt]
(q;q)_\infty(q^2;q^2)_\infty&  =\sum_{n=0}^{\infty}\sum_{\mid m\mid\leqslant[n/3]}
  (-1)^{n} q^{\frac{1}{2}(n^2-8m^2)+\frac{1}{2}n}.
  \end{align}

  For $K(z;q)$, Garvan \cite{Garvan-2014} established the following identity.
  \begin{align}
  &(zq^2,z^{-1}q^2,q^2;q^2)_{\infty}
  K(z;q) \nonumber \\[3pt] \nonumber
  &\hskip 0.2cm =(zq^2,z^{-1}q^2,q^2;q^2)_{\infty}
  \sum_{n=0}^{\infty}\frac{(-1)^nq^{n^2}(q;q^2)_n}
  {(zq^2;q^2)_n(z^{-1}q^2;q^2)_n}\\[3pt]
  &\hskip 0.2cm =\sum_{n=0}^{\infty}\Big(\sum_{m=0}^n(-1)^nz^{m-n}
  q^{\frac{1}{2}(2n^2-m^2)+\frac{1}{2}(2n-m)}
  +\sum_{m=1}^n(-1)^nz^{n-m+1}q^{\frac{1}{2}(2n^2-m^2)
  +\frac{1}{2}(2n+m)}\Big). \label{Garvan4}
\end{align}
 Setting $z=1$ in \eqref{Garvan4}, and noting that  $K(1;q)=(q;q)_\infty/(q^2;q^2)^2_\infty$,  we   obtain  the following   Hecke-Rogers type identity of Bressoud \cite{Bressoud-1986}.
\begin{equation}
(q;q)_\infty(q^2;q^2)_\infty=\sum_{n=0}^\infty \sum_{m=-n}^n(-1)^n q^{\frac{1}{2}(2n^2-m^2)  +\frac{1}{2}(2n+m)}.
\end{equation}

In \cite{Garvan-2014}, Garvan showed that \eqref{Garvan2} and \eqref{Garvan3} are equivalent with the aid of a transform of Mine \cite{Milne-1994}. He then proved \eqref{Garvan1}, \eqref{Garvan2} and \eqref{Garvan4} are valid   by showing the coefficients of $z^k$ on both sides of these   identities are equal.    At the end of his paper, Garvan proposed a problem of finding  direct proofs of these three identities by the method of  Bailey pairs. In this paper, we will show these  identities hold by applying Bailey's transform with a conjugate Bailey pair  of Warnaar \cite{Warnaar-2003} and three Bailey pairs derived from three identities  of Slater \cite{Slater-1949}. In particular, we obtain a more compact form of Hecke-Rogers identity for $R(z;q)$.

\begin{thm} \label{Garvan-n1-th}  We have
\begin{align}\label{Garvan-n1}\nonumber
&(zq, z^{-1}q,q;q)_{\infty}R(z;q)\\[3pt]
  &=\sum_{n=0}^{\infty}\left(\sum_{j=0}^{[n/2]}(-1)^{n+j}z^{n-3j}
 q^{\frac{1}{2}(n^2-3j^2)+\frac{1}{2}(n-j)}
  +\sum_{j=1}^{[n/2]}(-1)^{n+j}z^{n-3j+1}
 q^{\frac{1}{2}(n^2-3j^2)+\frac{1}{2}(n+j)}\right)\\[5pt]
 &=\sum_{n=0}^{\infty}\sum_{0\leq j<n/3}(-1)^{n+j}(z^{n-3j}+z^{3j-n})
 q^{\frac{1}{2}(n^2-3j^2)+\frac{1}{2}(n-j)}
 +\sum_{j=0}^{\infty}q^{3j^2+j}\nonumber\\[5pt]
 &\quad +\sum_{n=0}^{\infty}\sum_{1\leq j<(n+1)/3}(-1)^{n+j}(z^{n-3j+1}+z^{3j-n-1})
 q^{\frac{1}{2}(n^2-3j^2)+\frac{1}{2}(n+j)}  -\sum_{j=1}^{\infty}q^{3j^2-j}.  \label{Garvan-n1-a}
\end{align}
\end{thm}
Note that Garvan's identity \eqref{Garvan1}  can  immediately follow from   \eqref{Garvan-n1}. To be specific, if we replace $z$ by $z^{-1}$ in \eqref{Garvan-n1}, we get
 \begin{align}  \label{tempt}
&(zq, z^{-1}q,q;q)_{\infty}R(z^{-1};q)\\ \nonumber
  &=\sum_{n=0}^{\infty}\left(\sum_{j=0}^{[n/2]}(-1)^{n+j}z^{-n+3j}
 q^{\frac{1}{2}(n^2-3j^2)+\frac{1}{2}(n-j)}
  +\sum_{j=1}^{[n/2]}(-1)^{n+j}z^{-n+3j-1}
 q^{\frac{1}{2}(n^2-3j^2)+\frac{1}{2}(n+j)}\right).
\end{align}
  Combining \eqref{Garvan-n1} and \eqref{tempt},  and noting    $R(z;q)=R(z^{-1};q)$, we get    Garvan's identity \eqref{Garvan1}.

Using the Bailey machinery, we obtain the following Hecke-Rogers identity for $g_2(z;q)$.

\begin{thm} \label{Garvan-n2-th}  We have
\begin{equation}
 \begin{split}\label{Garvan-n2}
  &(1-z)(zq, z^{-1}q,q;q)_{\infty}g_2(z;q)\\
    &=\sum_{n=0}^{\infty}\left(\sum_{m=0}^{[n/2]}
    (-1)^{m+n}z^{n-2m}q^{\frac{1}{2}(n^2-2m^2)+\frac{1}{2}n}
   +\sum_{m=1}^{[n/2]}
    (-1)^{m+n}z^{2m-n-1}q^{\frac{1}{2}(n^2-2m^2)+\frac{1}{2}n}\right).
 \end{split}
 \end{equation}
\end{thm}
When  substitute   \eqref{Garvan-n2} into the relation \eqref{r-h-g2}, and employ Jacobi's triple product identity,  we derive  Garvan's identity \eqref{Garvan2}.

More generally, we generalized  these two-variable Hecke-Rogers identities to infinite families. The next theorem embeds  \eqref{Garvan-n1-a} into an infinite family.

\begin{thm}\label{Garvan-n1-gen}
For $k\geq 2$,  we have
    \begin{align}\label{gener1}
&(zq, z^{-1}q,q;q)_{\infty}\sum_{n=0}^{\infty}
\frac{(q;q)_{2n}q^n}{(zq;q)_n(z^{-1}q;q)_n}\sum_{n_1,\cdots,n_{k-1}=0}^{\infty}
\frac{q^{N_1^2+\cdots+N_{k-2}^2+2N_{k-1}^2+N_1+\cdots+N_{k-2}}}
{(q;q)_{n-N_1}(q;q)_{n_1}\cdots(q;q)_{n_{k-2}}(q;q)_{2n_{k-1}}}\nonumber \\[3pt]
&=\sum_{n=0}^{\infty}\sum_{0\leq j<n/3}
(-1)^{n+j}(z^{n-3j}+z^{3j-n})q^{\frac{1}{2}(n^2-3j^2)+
\frac{1}{2}(n-j)+3(k-1)j(3j+1)}\nonumber  \\
&\quad+\sum_{n=0}^{\infty}\sum_{1\leq j
<(n+1)/3}(-1)^{n+j}(z^{n-3j+1}+z^{3j-n-1})q^{\frac{1}{2}(n^2-3j^2)+
\frac{1}{2}(n+j)+3(k-1)j(3j-1)}\nonumber \\
&\quad +\sum_{j=0}^{\infty}q^{(3k-2)j(3j+1)}
-\sum_{j=1}^{\infty}q^{(3k-2)j(3j-1)},
    \end{align}   
  where  $N_j=n_j+n_{j+1}+\cdots+n_{k-1}$.
\end{thm}

We extends \eqref{Garvan-n2} into an infinite family.

\begin{thm}\label{gener2-thm}
For $k\geq 2$, we have 
    \begin{align}\label{gener2}
&(zq, z^{-1}q,q;q)_{\infty}\sum_{n=0}^{\infty}
\frac{(q^2;q)_{2n}q^n}{(zq;q)_n(z^{-1}q;q)_{n+1}}\nonumber \\[2pt]
&\hskip 3cm \times\sum_{n_1,\cdots,n_{k-1}=0}^{\infty}
\frac{q^{N_1^2+\cdots+N_{k-2}^2+\frac{3}{2}N_{k-1}^2
+2(N_1+\cdots+N_{k-2})
+\frac{3}{2}N_{k-1}}(-q;q)_{n_{k-1}}}
{(q;q)_{n-N_1}(q;q)_{n_1}\cdots(q;q)_{n_{k-2}}
(q^2;q)_{2n_{k-1}}}\nonumber\\[2pt]
&=\sum_{n=0}^{\infty}\left(\sum_{m=0}^{[n/2]}
(-1)^{m+n}z^{n-2m}q^{\frac{1}{2}(n^2-2m^2)+\frac{1}{2}n
+4(k-1)(m^2+m)}\right. \nonumber\\ 
&\quad\hskip1cm\left.+\sum_{m=1}^{[n/2]}(-1)^{m+n}z^{2m-n-1}
q^{\frac{1}{2}(n^2-2m^2)+\frac{1}{2}n+4(k-1)(m^2-m)}\right),
    \end{align} 
where   $N_j=n_j+n_{j+1}+\cdots+n_{k-1}$.
\end{thm}

Finally, we generalize \eqref{Garvan4} to an infinite family.

\begin{thm}\label{gener3-thm}
For $k\geq 2$, we have
\begin{align}\label{gener3}
&(zq^2, z^{-1}q^2,q^2;q^2)_{\infty}\sum_{n=0}^{\infty}
\frac{(q^2;q^2)_{2n}q^{2n}}{(zq^2;q^2)_n(z^{-1}q^2;q^2)_{n}}\nonumber \\[2pt]
&\hskip 4cm \times\sum_{n_1,\cdots,n_{k-1}=0}^{\infty}
\frac{(-1)^{n_{k-1}}q^{2N_1^2+\cdots+2N_{k-2}^2+3N_{k-1}^2
+2N_1+\cdots+2N_{k-2}}(q;q^2)_{n_{k-1}}}
{(q^2;q^2)_{n-N_1}(q^2;q^2)_{n_1}\cdots(q^2;q^2)_{n_{k-2}}
(q^2;q^2)_{2n_{k-1}}}\nonumber \\
&=\sum_{n=0}^{\infty}\left(\sum_{m=0}^n(-1)^{n}z^{m-n}
q^{\frac{1}{2}(2n^2-m^2)+\frac{1}{2}(2n-m)+2(k-1)(m^2+m)}\right. \nonumber \\
&\quad \hskip 1cm\left.+\sum_{m=1}^n(-1)^{n}z^{n-m+1}
q^{\frac{1}{2}(2n^2-m^2)+\frac{1}{2}(2n+m)+2(k-1)(m^2-m)}\right),
\end{align} 
where  $N_j=n_j+n_{j+1}+\cdots+n_{k-1}$.
\end{thm}

\section{Conjugate Bailey pairs and Bailey pairs  }

In this section, we   present some conjugate Bailey pairs and Bailey pairs which are needed in the proof of Hecke-Rogers identities for the universal mock theta functions. Recall that a pair of sequences $(\delta_n, \gamma_n)$ is called a conjugate Bailey pair relative to $a$  if they satisfy
\begin{equation}
\gamma_n=\sum_{r=n}^\infty\frac{\delta_r}{(q;q)_{r-n}(aq;q)_{r+n}}.
\end{equation}
A pair of sequences $(\alpha_n, \beta_n)$ is called a Bailey pair relative to $a$ if they satisfy
\begin{equation}
\beta_n=\sum_{r=0}^n\frac{\alpha_r}{(q;q)_{n-r}(aq;q)_{n+r}}.
\end{equation}

From the definitions of Bailey pairs and conjugate Bailey pairs, one could easily obtain the following transform \cite{Bailey-1949}.

\begin{thm}[The Bailey transform] If $(\alpha_n,\beta_n)$ is a Bailey pair relative to $a$ and $(\gamma_n,\delta_n)$ is a conjugate Bailey pair relative to $a$. Then we have
\begin{equation*}
    \sum_{n=0}^{\infty}\alpha_{n}\gamma_{n}
    =\sum_{n=0}^{\infty}\beta_{n}\delta_{n}.
\end{equation*}

\end{thm}

 We first state a conjugate Bailey pair which is needed in the proof of all of Garvan's three identities.  The proof of this conjugate Bailey pair is mainly based on the following generalization
  of Jacobi's triple product identity given by Warnaar \cite{Warnaar-2003}. Andrews and Warnaar \cite{Andrews-Warnaar-2007} provided an alternative proof of this identity.
 \begin{equation}\label{Jacobi}
    1+\sum_{r=1}^{\infty}(-1)^rq^{r\choose 2}(a^r+b^r)
    =(a,b,q;q)_{\infty}
    \sum_{r=0}^{\infty}\frac{(ab/q;q)_{2r}q^r}{(q,a,b,ab;q)_r}.
 \end{equation}
 Setting $a=z$ and $b=z^{-1}q$ in \eqref{Jacobi}, we obtain  Jacobi's  triple product \cite[p.15]{Gasper-2004}
\begin{equation}\label{Jacobi-o}
1+\sum_{r=1}^{\infty}(-1)^rq^{r\choose 2}(z^r+z^{-r}q^r)
=(z,z^{-1}q,q;q)_{\infty}.
\end{equation}

 \begin{lem}\label{lem1}The pair of sequences $(\delta_n, \gamma_n)$ forms a conjugate Bailey pair relative to $ab$, where
  \begin{equation*}
 \delta_n=(aq,b,q;q)_{\infty}
   \frac{(ab;q)_{2n}q^n}{(aq;q)_{n}(b;q)_n},
 \end{equation*}
 and
 \begin{equation*}
    \gamma_n=\frac{(1-ab) q^{n}}{(1-abq^{2n})}\left(1+\sum_{r=1}^{\infty}
    (-1)^rq^{r \choose 2}((aq^{n+1})^r+(bq^n)^r)\right).
 \end{equation*}
 \end{lem}

 \pf
 By the definition of conjugate Bailey pair, we see that
 \begin{align*}
    \gamma_n&=\sum_{r=n}^{\infty}\frac{\delta_r}{(q;q)_{r-n}(abq;q)_{r+n}}\\[3pt]
    &=(aq,b,q;q)_{\infty}\sum_{r=n}^{\infty}
    \frac{(ab;q)_{2r}q^r}{(q;q)_{r-n}(abq;q)_{r+n}(aq;q)_{r}(b;q)_r}\\[3pt]
    &=(aq,b,q;q)_{\infty}\sum_{r=0}^{\infty}
    \frac{(ab;q)_{2n+2r}q^{r+n}}{(q;q)_{r}(abq;q)_{r+2n}(aq;q)_{r+n}(b;q)_{r+n}}\\[3pt]
    &=(aq,b,q;q)_{\infty}\frac{(ab;q)_{2n}q^n}
    {(abq;q)_{2n}(aq;q)_{n}(b;q)_{n}}\sum_{r=0}^{\infty}
    \frac{(abq^{2n};q)_{2r}q^{r}}{(q,abq^{2n+1},aq^{n+1},bq^n;q)_{r}}\\[3pt]
    &=\frac{(1-ab)q^n}{1-abq^{2n}}(aq^{n+1},bq^n,q;q)_\infty\sum_{r=0}^{\infty}
    \frac{(abq^{2n};q)_{2r}q^{r}}{(q,abq^{2n+1},aq^{n+1},bq^n;q)_{r}}.
    \end{align*}
 Using \eqref{Jacobi} with   $a$ replaced by $aq^{n+1}$ and $b$ by $bq^n$ in the summation of the above identity, we derive that 
    \begin{align*}
     \gamma_n&=\frac{(1-ab)q^n}{1-abq^{2n}}(aq^{n+1},bq^n,q;q)_\infty \frac{1}{(aq^{n+1},bq^n,q;q)_\infty} \\[3pt]
    &\hskip 2cm \times \left(1+\sum_{r=1}^{\infty}
    (-1)^rq^{r \choose 2}((aq^{n+1})^r+(bq^n)^r)\right)\\[3pt]
    &=\frac{(1-ab)q^n}{1-abq^{2n}}\left(1+\sum_{r=1}^{\infty}
    (-1)^rq^{r \choose 2}((aq^{n+1})^r+(bq^n)^r)\right),
 \end{align*}
 as desired. This completes the proof.  \qed

To derive \eqref{Garvan-n1-a}, we also need  the following Bailey pair given by Warnnar \cite[p.375]{Warnaar-2003}, which can be established by a identity of Slater \cite[Eq.(3.4)]{Slater-1949}.

\begin{lem}\label{Garvan-BL-1}The pair of sequences $(\alpha_n,\beta_n)$ forms a Bailey pair relative to $q$, where
\begin{equation*}
\begin{split}
 \alpha_{3n}&=q^{(3n-2)n}\frac{1-q^{6n+1}}{1-q},  \quad \alpha_{3n+1}=0,
 \quad \alpha_{3n+2}=-q^{(3n+2)n}\frac{1-q^{6n+5}}{1-q}, \\[3pt]
 \beta_n&=\frac{q^{n(n-1)}}{(q;q)_{2n}}.
 \end{split}
\end{equation*}
\end{lem}

 The following Bailey pair is necessary in the proof of \eqref{Garvan-n2}. The proof of this Bailey pair is mainly based on a special case of $_6\psi_6$ given by Slater \cite[Eq.(4.2)]{Slater-1949}:
\begin{align}\label{Slater}
  &\sum_{r=-[n/2]}^{[n/2]}\frac{(1-aq^{4r})(q^{-n};q)_{2r}
  (d;q^2)_r(e;q^2)_r}{(1-a)(aq^{n+1};q)_{2r}(aq^2/d;q^2)_r(aq^2/e;q^2)_r}
 \left( \frac{a^{2}q^{2n+1}}{de}\right)^r \\\nonumber
  &\hskip1cm=\frac{(q^2/a,aq/d,aq/e,aq^2/de;q^2)_{\infty}}
  {(q,q^2/d,q^2/e,a^2q/de;q^2)_\infty}\cdot
  \frac{(q;q)_n(aq;q)_n(a^2q/de;q^2)_n}{(aq;q^2)_n(aq/d;q)_n(aq/e;q)_n}.
\end{align}

 \begin{lem}\label{Garvan-BL-2}
 The pair of sequences $(\alpha_n,\beta_n)$ forms a Bailey pair relative to $q^2$, where
 \begin{equation}\label{Bailey pair1}
    \begin{split}
     \alpha_{2n}&=(-1)^nq^{n^2-n}\frac{1-q^{4n+2}}{1-q^2},
     \ \ \ \alpha_{2n+1}=0,\\
     \beta_n&=\frac{(-q;q)_nq^{n\choose 2}}{(q^2;q)_{2n}}.
    \end{split}
 \end{equation}
 \end{lem}

\pf In \eqref{Slater}, let  $a=d=q^2$,  and
$e\rightarrow0$, noting  that $1/(q^2;q^2)_n=0$ when  $n<0$, we obtain
\[\sum_{r=0}^{[n/2]}\frac{(1-q^{4r+2})(q^{-n};q)_{2r}(-1)^rq^{2nr-r^2}}
{(1-q^2)(q^{n+3};q)_{2r}}=\frac{(q^3;q)_nq^{n\choose 2}}{(q^3;q^2)_n}.\]
Divide   both sides of the above identity  by $(q;q)_n(q^3;q)_n$   to get
\begin{equation*}
    \sum_{r=0}^{[n/2]}\frac{1}{(q;q)_{n-2r}(q^3;q)_{n+2r}}
    \frac{(1-q^{4r+2})(-1)^rq^{r^2-r}}{1-q^2}=\frac{(-q;q)_nq^{n\choose 2}}{(q^2;q)_{2n}},
\end{equation*}
which yields the desired Bailey pair. This completes the proof. \qed

To verify \eqref{Garvan4}, we also require the following Bailey pair.

 \begin{lem}\label{Garvan-BL-3}
 The pair of sequences $(\alpha_n,\beta_n)$ forms a Bailey pair relative to $q$ , where
 \begin{equation}\label{Bailey pair2}
    \begin{split}
      \alpha_n&=(-1)^n q^{\frac{n^2-3n}{4}}\frac{1-q^{2n+1}}{1-q},\\[3pt]
      \beta_n&=\frac{(-1)^nq^{\frac{n^2}{2}-n}(q^{\frac{1}{2}};q)_n}{(q;q)_{2n}}.
    \end{split}
 \end{equation}
 \end{lem}

 \pf In \eqref{Slater}, set $a=q,\,d=-q^{3/2}$
  and $e\rightarrow 0$ to get
  \begin{align*}
  \sum_{r=-[n/2]}^{[n/2]}\frac{(1-q^{4r+1})(q^{-n};q)_{2r}
  q^{2nr-\frac{1}{2}r-r^2}}
    {(1-q)(q^{n+2};q)_{2r}}=\frac{(q;q)_n(q^2;q)_n(-1)^nq^{\frac{n^2}{2}-n}}
    {(q^2;q^2)_n(-q^{\frac{1}{2}};q)_n}.
  \end{align*}
  Dividing  both sides of the above identity  by $(q;q)_n(q^2;q)_n$,  we obtain
  \begin{equation*}
    \sum_{r=-[n/2]}^{[n/2]}\frac{(1-q^{4r+1})q^{r^2-\frac{3}{2}r}}
    {(q;q)_{n+2r+1}(q;q)_{n-2r}}
    =\frac{(-1)^nq^{\frac{n^2}{2}-n}(q^{\frac{1}{2}};q)_n}{(q;q)_{2n}}.
  \end{equation*}
Simplifying the left hand side of the above identity yields
\begin{align*}
 &\sum_{r=-[n/2]}^{[n/2]}\frac{(1-q^{4r+1})q^{r^2-\frac{3}{2}r}}
    {(q;q)_{n+2r+1}(q;q)_{n-2r}}\\[5pt]
&\hskip 0.2cm= \frac{1}{(q^2;q)_n(q;q)_n}+\sum_{r=1}^{[n/2]}\frac{(1-q^{4r+1})
q^{r^2-\frac{3}{2}r}}
    {(q;q)_{n+2r+1}(q;q)_{n-2r}}+\sum_{r=-[n/2]}^{-1}\frac{(1-q^{4r+1})
    q^{r^2-\frac{3}{2}r}}{(q;q)_{n+2r+1}(q;q)_{n-2r}}\\[5pt]
    &=\frac{1}{(q^2;q)_n(q;q)_n}+\sum_{r=1}^{[n/2]}
    \frac{(1-q^{4r+1})q^{r^2-\frac{3}{2}r}}
    {(q;q)_{n+2r+1}(q;q)_{n-2r}}-\sum_{r=1}^{[n/2]}
    \frac{(1-q^{4r-1})q^{r^2-\frac{5}{2}r+1}}
    {(q;q)_{n-2r+1}(q;q)_{n+2r}}\\[5pt]
        &=\sum_{r=0}^n\frac{1}{(q;q)_{n-r}(q^2;q)_{n+r}}
    \frac{(1-q^{2r+1})(-1)^rq^{\frac{r^2-3r}{4}}}{1-q}.
\end{align*}
Hence, we have 
\[\sum_{r=0}^n\frac{1}{(q;q)_{n-r}(q^2;q)_{n+r}}
    \frac{(1-q^{2r+1})(-1)^rq^{\frac{r^2-3r}{4}}}{1-q}=\frac{(-1)^nq^{\frac{n^2}{2}-n}(q^{\frac{1}{2}};q)_n}{(q;q)_{2n}},\]
which implies the desired Bailey pair. This completes the proof. \qed

To get an infinite family of Hecke-Rogers identities,
  we need the following limiting form of Bailey's lemma \cite{Andrews-1984-a, Paule-1985}.

\begin{thm}\label{BL-c-1}
Supposed   $(\alpha_n,\beta_n)$ is a Bailey pair relative to $a$. Then $(\alpha'_n,\beta'_n)$ is a new Bailey pair relative to $a$, where

\begin{equation*}\label{a}
\alpha_n' =a^n q^{n^2}\alpha_n
\end{equation*}
and
\begin{align*}\label{b}
\beta_n' &=\sum_{j=0}^{\infty}\frac{a^jq^{j^2}}
{(q;q)_{n-j}} \beta_j .
\end{align*}
\end{thm}
Iterating this lemma leads to the Bailey chain.
We now produce the following infinite sequence of Bailey pairs by iterating the Bailey pair in Lemma \ref{Garvan-BL-1} along this Bailey chain.

\begin{lem}\label{Garvan-BL-1-gen} For $k\geq 2$,   $(\alpha^{(k)}_n, \beta_n^{(k)}
)$ form   Bailey pairs relative to $q$, where
\begin{equation*}
  \begin{split}
    \alpha_{3n}^{(k)}&=q^{3(3k-2)n^2+(3k-5)n}
    \frac{(1-q^{6n+1})}{(1-q)}, \ \ \alpha_{3n+1}^{(k)}=0,\\[3pt]
     \alpha_{3n+2}^{(k)}&=-q^{3(3k-2)n^2+(15k-13)n+6(k-1)}
    \frac{(1-q^{6n+5})}{(1-q)},\\[3pt]
    \beta_n^{(k)}&=\sum_{n_1,\cdots,n_{k-1}=0}^{\infty}
\frac{q^{N_1^2+\cdots+N_{k-2}^2+2N_{k-1}^2+N_1+\cdots+N_{k-2}}}
{(q;q)_{n-N_1}(q;q)_{n_1}\cdots(q;q)_{n_{k-2}}(q;q)_{2n_{k-1}}},
  \end{split}
\end{equation*}
where  $N_j=n_j+n_{j+1}+\cdots+n_{k-1}$.
\end{lem}
\pf We proceed by induction on $k$. When $k=2$. Applying    Theorem  \ref{BL-c-1} to the Bailey pair in Lemma \ref{Garvan-BL-1},  we obtain a new Bailey pair $(\alpha^{(2)}_n,\beta^{(2)}_n)$  relative to $q$, where
\begin{align*}
    \alpha_{3n}^{(2)}&=q^{9n^2+3n}\cdot q^{(3n-2)n}\frac{1-q^{6n+1}}{1-q}
    =q^{12n^2+n}\frac{1-q^{6n+1}}{1-q}, \quad 
      \alpha_{3n+1}^{(2)}=0,\\[5pt]
    \alpha_{3n+2}^{(2)}&=q^{(3n+2)(3n+3)}\cdot \left(-q^{(3n+2)n}\frac{1-q^{6n+5}}{1-q}\right)
    =-q^{12n^2+17n+6}\frac{1-q^{6n+5}}{1-q}
\end{align*}
and
\begin{align*}
    \beta_n^{(2)}=\sum_{j=0}^n\frac{q^{2j^2}}{(q;q)_{n-j}(q;q)_{2n}}
    =\sum_{j=0}^\infty\frac{q^{2j^2}}{(q;q)_{n-j}(q;q)_{2n}},
\end{align*}
as desired.

Assume that we have a Bailey pair $(\alpha^{(k-1)}_n,\beta^{(k-1)}_n)$ relative to $q$, where
\begin{equation*}
  \begin{split}
    \alpha_{3n}^{(k-1)}&=q^{3(3k-5)n^2+(3k-8)n}
    \frac{(1-q^{6n+1})}{(1-q)},  \ \ \alpha_{3n+1}^{(k-1)}=0,\\[3pt]
     \alpha_{3n+2}^{(k-1)}&=-q^{3(3k-5)n^2+(15k-28)n+6(k-2)}
\frac{(1-q^{6n+5})}{(1-q)},\\[3pt]
    \beta_n^{(k-1)}&=\sum_{n_1',\cdots,n_{k-2}'=0}^{\infty}
\frac{q^{N_1'^2+\cdots+N_{k-3}'^2+2N_{k-2}'^2+N_1'+\cdots+N_{k-3}'}}
{(q;q)_{n-N_1'}(q;q)_{n_1'}\cdots(q;q)_{n_{k-3}'}(q;q)_{2n_{k-2}'}},
  \end{split} 
\end{equation*}
 where $N_j'=n_j'+n_{j+1}'+\cdots+n_{k-2}'$. 

When we apply Theorem \ref{BL-c-1} to this Bailey pair, we  arrive at a new Bailey pair $(\alpha^{(k)}_n,\beta^{(k)}_n)$ relative to $q$, where for $\alpha^{(k)}_n$, we have
\begin{align*}
\alpha_{3n}^{(k)}&=q^{9n^2+3n}\alpha_{3n}^{(k-1)}\\[3pt]
&=q^{9n^2+3n+3(3k-5)n^2+(3k-8)n}\frac{(1-q^{6n+1})}{(1-q)}\\[3pt]
&=q^{3(3k-2)n^2+(3k-5)n}\frac{(1-q^{6n+1})}{(1-q)},\\[5pt]
\alpha_{3n+1}^{(k)}&=0, \\[5pt]
\alpha_{3n+2}^{(k)}&=q^{3(n+1)(3n+2)}\alpha_{3n+2}^{(k-1)}\\[3pt]
&=-q^{3(n+1)(3n+2)+3(3k-5)n^2+(15k-28)n+6(k-2)}
\frac{(1-q^{6n+5})}{(1-q)}\\[3pt]
&=-q^{3(3k-2)n^2+(15k-13)n+6(k-1)}\frac{(1-q^{6n+5})}{(1-q)}.
\end{align*}
For $\beta^{(k)}_n$, we have
\begin{align*}
    \beta_{n}^{(k)}&=\sum_{n_{k-1}'=0}^\infty\frac{q^{n_{k-1}'^2+n_{k-1}'}}
    {(q;q)_{n-n_{k-1}'}}\beta_{n'_{k-1}}^{(k-1)}\\
    &=\sum_{n_{k-1}'=0}^\infty\frac{q^{n_{k-1}'^2+n_{k-1}'}}{(q;q)_{n-n_{k-1}'}}
    \sum_{n_1',\cdots,n_{k-2}'=0}^{\infty}
\frac{q^{N_1'^2+\cdots+N_{k-3}'^2+2N_{k-2}'^2+N_1'+\cdots+N_{k-3}'}}
{(q;q)_{n_{k-1}'-N_1'}(q;q)_{n_1'}\cdots(q;q)_{n_{k-3}'}
(q;q)_{2n_{k-2}'}}.
\end{align*}
Recall that $N_j'=n_j'+n_{j+1}'+\cdots+n_{k-2}'$ for $1\leq j\leq k-2$. Let
\[n_1=n_{k-1}'-N_1',\  n_2=n_1',\  n_3=n_2',\ \ldots,\  n_{k-1}=n_{k-2}',\] and set $N_j=n_j+n_{j+1}+\cdots+n_{k-1}$ for $1\leq j\leq k-1$. Obviously, we have
\[N_2=N_1',\ N_3=N_2',\ \ldots,\ N_{k-1}=N_{k-2}',\ N_1=n_{k-1}',\]
and so
 \begin{align*}
   \beta_{n}^{(k)}&=\sum_{n_1,\cdots,n_{k-1}=0}^{\infty}
\frac{q^{N_1^2+\cdots+N_{k-2}^2+2N_{k-1}^2+N_1+\cdots+N_{k-2}}}
{(q;q)_{n-N_1}(q;q)_{n_1}\cdots(q;q)_{n_{k-2}}(q;q)_{2n_{k-1}}}.
 \end{align*}
This completes the proof. \qed

The following infinite sequence of Bailey pairs can be obtained from the Bailey pair in Lemma \ref{Garvan-BL-2} along the same line.

\begin{lem}\label{Garvan-BL-2-gen}For $k\geq 2$,   $(\alpha^{(k)}_n, \beta_n^{(k)}
)$ form   Bailey pairs relative to $q^2$, where
\begin{align*}
      \alpha_{2n}^{(k)}&=(-1)^nq^{(4k-3)n^2+(4k-5)n}\frac{1-q^{4n+2}}{1-q^2},
 \quad  \alpha_{2n+1}^{(k)}=0,\\[5pt]
       \beta_{n}^{(k)}&=\sum_{n_1,\cdots,n_{k-1}=0}^{\infty}
\frac{q^{N_1^2+\cdots+N_{k-2}^2+\frac{3}{2}N_{k-1}^2+2(N_1+\cdots+N_{k-2})
+\frac{3}{2}N_{k-1}}(-q;q)_{n_{k-1}}}
{(q;q)_{n-N_1}(q;q)_{n_1}\cdots(q;q)_{n_{k-2}}(q^2;q)_{2n_{k-1}}}.
\end{align*}
\end{lem}

Similarly, we obtain the following infinite sequence of Bailey pairs from the Bailey pair in Lemma \ref{Garvan-BL-3} which are used to establish the infinite family in Theorem \ref{gener3-thm}.

\begin{lem}\label{Garvan-BL-3-gen}For $k\geq 2$,   $(\alpha^{(k)}_n, \beta_n^{(k)}
)$ form   Bailey pairs relative to $q$, where
\begin{align*}
      \alpha_n^{(k)}&=(-1)^nq^{\frac{4k-3}{4}n^2+\frac{4k-7}{4}n}
   \frac{ 1-q^{2n+1}}{1-q},\\[5pt]
    \beta_n^{(k)}&=\sum_{n_1,\cdots,n_{k-1}=0}^{\infty}
\frac{(-1)^{n_{k-1}}q^{N_1^2+\cdots+N_{k-2}^2+\frac{3}{2}N_{k-1}^2
+N_1+\cdots+N_{k-2}}(q^{\frac{1}{2}};q)_{n_{k-1}}}
{(q;q)_{n-N_1}(q;q)_{n_1}\cdots(q;q)_{n_{k-2}}
(q;q)_{2n_{k-1}}}.
\end{align*}
\end{lem}

\section{A two-variable Hecke-Rogers identity for $R(z;q)$ and   generalization}

In this section, we first give a proof of \eqref{Garvan-n1-a} in Theorem \ref{Garvan-n1-th} by using the Bailey transform with conjugate Bailey pair in Lemma \ref{lem1} and Bailey pair in Lemma \ref{Garvan-BL-1}. We then show that \eqref{Garvan-n1} can be derived from \eqref{Garvan-n1-a} with the aid of  two identities given by Garvan \cite{Garvan-2014}. Finally, we give a proof of Theorem \ref{Garvan-n1-gen}.

\medskip

 \noindent{\it Proof of Theorem \ref{Garvan-n1-th}, Eq. \eqref{Garvan-n1-a}.}  Combining the Bailey
 transform with the Bailey pairs in Lemma \ref{Garvan-BL-1} and the conjugate Bailey pair in Lemma \ref{lem1} with $a=z$, and $b=z^{-1}q$, we obtain
\begin{align}
\sum_{n=0}^{\infty}\beta_n\delta_n&=(zq,z^{-1}q,q;q)_{\infty}
\sum_{n=0}^{\infty}\frac{q^{n^2}}{(zq;q)_n(z^{-1}q;q)_n}\nonumber\\
&=\sum_{j=0}^{\infty}\alpha_{3j}\gamma_{3j}+\sum_{j=1}^{\infty}\alpha_{3j-1}\gamma_{3j-1}\nonumber\\
&=\sum_{j=0}^{\infty}q^{j(3j-2)+3j}\left(1+\sum_{n=1}^{\infty}(-1)^nq^{n\choose 2}
((zq^{3j+1})^n+(z^{-1}q^{3j+1})^n)\right)\nonumber\\
&\quad -\sum_{j=1}^{\infty}q^{(j-1)(3j-1)+(3j-1)}\left(1+\sum_{n=1}^{\infty}(-1)^nq^{n\choose 2}
((zq^{3j})^n+(z^{-1}q^{3j})^n)\right)\nonumber
\end{align}
\begin{align}
&=\sum_{j=0}^{\infty}\sum_{n=1}^{\infty}(-1)^n(z^n+z^{-n})
 q^{{n+1\choose 2}+3j^2+j+3nj}+\sum_{j=0}^{\infty}q^{3j^2+j}\nonumber\\
 &\quad + \sum_{j=1}^{\infty}\sum_{n=1}^{\infty}(-1)^{n-1}(z^n+z^{-n})
 q^{{n\choose 2}+3j^2-j+3nj}-\sum_{j=1}^{\infty}q^{3j^2-j}\label{garvan}.
 \end{align}
Replace  $n$ by $n-3j$ in the first sum  of \eqref{garvan} and interchange  the order of  summation to get
\begin{align}\label{sum1}
&\sum_{j=0}^{\infty}\sum_{n=1}^{\infty}(-1)^n(z^n+z^{-n})
 q^{{n+1\choose 2}+3j^2+j+3nj}\nonumber \\[2pt]
 &\quad =\sum_{n=0}^{\infty}\sum_{0\leq j<n/3}(-1)^{n+j}(z^{n-3j}+z^{3j-n})
 q^{\frac{1}{2}(n^2-3j^2)+\frac{1}{2}(n-j)}.
 \end{align}
Replace $n$ by $n-3j+1$ in
 the third sum of \eqref{garvan} and interchange the order of  summation to yield
 \begin{align}\label{sum2}
 &\sum_{j=1}^{\infty}\sum_{n=1}^{\infty}(-1)^{n-1}(z^n+z^{-n})
 q^{{n\choose 2}+3j^2-j+3nj}\nonumber\\[2pt]
 &=\sum_{n=0}^{\infty}\sum_{1\leq j<(n+1)/3}(-1)^{n+j}(z^{n-3j+1}+z^{3j-n-1})
 q^{\frac{1}{2}(n^2-3j^2)+\frac{1}{2}(n+j)}.
 \end{align}
 Plugging  \eqref{sum1} and \eqref{sum2} into \eqref{garvan}, we get the desired identity \eqref{Garvan-n1-a}.\qed

We are now in position to show \eqref{Garvan-n1} in Theorem \ref{Garvan-n1-th} holds.

 \noindent{\it Proof of Theorem \ref{Garvan-n1-th}, Eq. \eqref{Garvan-n1}.}  Obviously, it suffices to show that
 \begin{align}
&\sum_{n=0}^{\infty}\sum_{0\leq j<n/3}(-1)^{n+j}(z^{n-3j}+z^{3j-n})
 q^{\frac{1}{2}(n^2-3j^2)+\frac{1}{2}(n-j)}+\sum_{j=0}^{\infty}q^{3j^2+j}\label{Garvan-n1-pf}\\
 &\quad +\sum_{n=0}^{\infty}\sum_{1\leq j<(n+1)/3}(-1)^{n+j}(z^{n-3j+1}+z^{3j-n-1})
 q^{\frac{1}{2}(n^2-3j^2)+\frac{1}{2}(n+j)} -\sum_{j=1}^{\infty}q^{3j^2-j} \nonumber\\[5pt]
&=\sum_{n=0}^{\infty}\left(\sum_{j=0}^{[n/2]}(-1)^{n+j}z^{n-3j}
 q^{\frac{1}{2}(n^2-3j^2)+\frac{1}{2}(n-j)}
  +\sum_{j=1}^{[n/2]}(-1)^{n+j}z^{n-3j+1}
 q^{\frac{1}{2}(n^2-3j^2)+\frac{1}{2}(n+j)}\right). \nonumber
 \end{align}

 To this end, we first show that
 \begin{align}\label{p1}
 &\sum_{n=0}^{\infty}\sum_{0\leq j<n/3}(-1)^{n+j}(z^{n-3j}+z^{3j-n})
 q^{\frac{1}{2}(n^2-3j^2)+\frac{1}{2}(n-j)} +\sum_{j=0}^{\infty}q^{3j^2+j} \nonumber \\[3pt]
 & \quad =\sum_{n=0}^{\infty}\sum_{j=0}^{[n/2]}(-1)^{n+j}z^{n-3j}
 q^{\frac{1}{2}(n^2-3j^2)+\frac{1}{2}(n-j)}.
 \end{align}

Using  the following identity given by Garvan \cite[Eq. (2.14)]{Garvan-2014},
\begin{align} \label{gposi-1}
&\sum_{n=0}^{\infty}\sum_{0\leq j<n/3}(-1)^{n+j} z^{n-3j}
   q^{\frac{1}{2}(n^2-3j^2)+\frac{1}{2}(n-j)}\nonumber\\[3pt]
  & =\sum_{n=0}^{\infty}\sum_{n/3< j\leq n/2}(-1)^{n+j} z^{-n+3j}
   q^{\frac{1}{2}(n^2-3j^2)+\frac{1}{2}(n-j)},
\end{align}
we get
\begin{align*}
 &\sum_{n=0}^{\infty}\sum_{0\leq j<n/3}(-1)^{n+j}(z^{n-3j}+z^{3j-n})
 q^{\frac{1}{2}(n^2-3j^2)+\frac{1}{2}(n-j)} +\sum_{j=0}^{\infty}q^{3j^2+j}\\[3pt]
 &=\sum_{n=0}^{\infty}\sum_{0\leq j<n/3}(-1)^{n+j}z^{n-3j}
 q^{\frac{1}{2}(n^2-3j^2)+\frac{1}{2}(n-j)}\\
 &\quad +\sum_{n=0}^{\infty}\sum_{n/3< j\leq n/2}(-1)^{n+j}z^{n-3j}
 q^{\frac{1}{2}(n^2-3j^2)+\frac{1}{2}(n-j)}+\sum_{j=0}^{\infty}q^{3j^2+j}\\
 &=\sum_{n=0}^{\infty}\sum_{j=0}^{[n/2]}(-1)^{n+j}z^{n-3j}
 q^{\frac{1}{2}(n^2-3j^2)+\frac{1}{2}(n-j)},
\end{align*}
as desired.

Similarly, applying the following identity of Garvan \cite[Eq. (2.15)]{Garvan-2014}
\begin{align*}
    &\sum_{n=0}^{\infty}\sum_{1\leq j<(n+1)/3}(-1)^{n+j} z^{n-3j+1}
q^{\frac{1}{2}(n^2-3j^2)+\frac{1}{2}(n+j)}\nonumber\\[3pt]
  & \hskip 1cm  =\sum_{n=0}^{\infty}\sum_{(n+1)/3<j\leq n/2}(-1)^{n+j} z^{-n+3j-1}
  q^{\frac{1}{2}(n^2-3j^2)+\frac{1}{2}(n+j)},
\end{align*}
we can show that
\begin{align}\label{p2}
   &\sum_{n=0}^{\infty}\sum_{0\leq j<n+1/3}(-1)^{n+j}(z^{n-3j+1}+z^{3j-n-1})
 q^{\frac{1}{2}(n^2-3j^2)+\frac{1}{2}(n+j)}
 -\sum_{j=1}^{\infty}q^{3j^2-j}\nonumber\\
 &=\sum_{n=0}^\infty\sum_{j=1}^{[n/2]}(-1)^{n+j}z^{n-3j+1}
 q^{\frac{1}{2}(n^2-3j^2)+\frac{1}{2}(n+j)}.
\end{align}
Combining \eqref{p1} and \eqref{p2},  we obtain
\eqref{Garvan-n1-pf}. Thus we show \eqref{Garvan-n1} holds.\qed

For the reminder of this section, we show \eqref{gener1} holds  by using the Bailey pair technology. \medskip

\noindent{\it Proof of Theorem \ref{Garvan-n1-gen}.} Applying the Bailey transform with the conjugate Bailey pair in Lemma \ref{lem1} by setting $a=z$ and $b=z^{-1}q$ and the Bailey pair in Lemma \ref{Garvan-BL-1-gen},  we get
\begin{align}
    \sum_{n=0}^{\infty}\beta_n^{(k)}\delta_n&=
   (zq,z^{-1}q,q;q)_{\infty}\sum_{n=0}^{\infty}
\frac{(q;q)_{2n}q^n}{(zq;q)_n(z^{-1}q;q)_n}\nonumber \\[2pt]
&\hskip 4cm \times \sum_{n_1,\cdots,n_{k-1}=0}^{\infty}
\frac{q^{N_1^2+\cdots+N_{k-2}^2+2N_{k-1}^2+N_1+\cdots+N_{k-2}}}
{(q;q)_{n-N_1}(q;q)_{n_1}\cdots(q;q)_{n_{k-2}}(q;q)_{2n_{k-1}}}\nonumber\\
&=\sum_{j=0}^{\infty}\alpha_{3j}^{(k)}\gamma_{3j}
+\sum_{j=1}^{\infty}\alpha_{3j-1}^{(k)}\gamma_{3j-1}\nonumber\\
&=\sum_{j=0}^{\infty}q^{(3k-2)j(3j+1)}
\left(1+\sum_{n=1}^{\infty}(-1)^nq^{n+1\choose 2}((zq^{3j})^n+(z^{-1}q^{3j})^n)\right)\nonumber\\
&\quad-\sum_{j=1}^{\infty}q^{(3k-2)j(3j-1)}
\left(1+\sum_{n=1}^{\infty}(-1)^nq^{n\choose 2}((zq^{3j})^n+(z^{-1}q^{3j})^n)\right)\nonumber\\
&=\sum_{j=0}^{\infty}q^{(3k-2)j(3j+1)}+\sum_{j=0}^{\infty}\sum_{n=1}^{\infty}
(-1)^n(z^n+z^{-n})q^{{n+1\choose2}+3jn+(3k-2)j(3j+1)}\nonumber\\\label{garvan2}
&\quad-\sum_{j=1}^{\infty}q^{(3k-2)j(3j-1)}+\sum_{j=1}^{\infty}\sum_{n=1}^{\infty}
(-1)^{n-1}(z^n+z^{-n})q^{{n\choose2}+3jn+(3k-2)j(3j-1)}.
\end{align}
Replace  $n$ by $n-3j$ in the second sum  of \eqref{garvan2} and reverse   the order of summation to get
\begin{align}\label{sum3}
&\sum_{j=0}^{\infty}\sum_{n=1}^{\infty}(-1)^n(z^n+z^{-n})
 q^{{n+1\choose 2}+3nj+(3k-2)(3j^2+j)}\nonumber \\
 &=\sum_{n=0}^{\infty}\sum_{0\leq j<n/3}(-1)^{n+j}(z^{n-3j}+z^{3j-n})
 q^{\frac{1}{2}(n^2-3j^2)+\frac{1}{2}(n-j)+3(k-1)j(3j+1)}.
 \end{align}
Replace $n$ by $n-3j+1$ in
 the fourth sum of \eqref{garvan2} and interchange the order of summation. This gives
 \begin{align}\label{sum4}
 &\sum_{j=1}^{\infty}\sum_{n=1}^{\infty}(-1)^{n-1}(z^n+z^{-n})
 q^{{n\choose 2}+3nj+(3k-2)(3j^2-j)}\nonumber\\[2pt]
 &=\sum_{n=0}^{\infty}\sum_{1\leq j<(n+1)/3}(-1)^{n+j}(z^{n-3j+1}+z^{3j-n-1})
 q^{\frac{1}{2}(n^2-3j^2)+\frac{1}{2}(n+j)+3(k-1)j(3j-1)}.
 \end{align}
 Substituting \eqref{sum3} and \eqref{sum4} in \eqref{garvan2}, we are led to the desired identity \eqref{gener1}.\qed

\section{A two-variable Hecke-Rogers identity for $H(z;q)$ and generalization}

In this section  we first use the Bailey transform to show Theorem \ref{Garvan-n2-th}. We then derive Garvan's identity   \eqref{Garvan2} from  Theorem \ref{Garvan-n2-th}. Last, we give a proof of  Theorem \ref{gener2-thm}.

\medskip

\noindent{\it Proof of Theorem \ref{Garvan-n2-th}.}
Apply the Bailey transform with the conjugate Bailey pair in Lemma \ref{lem1} by setting $a=z^{-1}q$ and $b=zq$ and the Bailey pair in Lemma \ref{Garvan-BL-2} to yield
\begin{align}
\sum_{n\geq 0}\beta_n\delta_n &=(zq, z^{-1}q,q;q)_{\infty}
 \sum_{n=0}^{\infty}\frac{(-q;q)_nq^{n+1\choose 2}}{(z^{-1}q;q)_{n+1}(zq;q)_{n}}\nonumber \\
 &=(1-z)(zq, z^{-1}q,q;q)_{\infty}
 g_2(z;q)\nonumber \\
 &=\sum_{m=0}^{\infty}\alpha_{2m}\gamma_{2m} \nonumber\\
 &=\sum_{m=0}^{\infty}(-1)^mq^{m^2+m}\left(1+\sum_{n=1}^{\infty}
 (-1)^nq^{n+1\choose 2}((z^{-1}q^{2m+1})^n+(zq^{2m})^n)\right) \nonumber\\
 &=\sum_{m=0}^{\infty}\sum_{n=0}^{\infty}
(-1)^{m+n}z^{n}q^{{n+1\choose 2}+2mn+m^2+m}\nonumber\\\label{h1}
&\quad\ \ +\sum_{m=0}^{\infty}
\sum_{n=1}^{\infty}(-1)^{m+n}z^{-n}q^{{n+1\choose 2}+(2m+1)n+m^2+m}.
\end{align}
Shift  $n$ to $n-2m$ in the first sum of \eqref{h1} and interchange the order of summation. This gives
\begin{equation}\label{sum-h1}
  \sum_{m=0}^{\infty}\sum_{n=0}^{\infty}
(-1)^{m+n}z^{n}q^{{n+1\choose 2}+2mn+m^2+m}
= \sum_{n=0}^{\infty}\sum_{m=0}^{[n/2]}
    (-1)^{m+n}z^{n-2m}q^{\frac{1}{2}(n^2-2m^2)+\frac{1}{2}n}.
\end{equation}
Shift  $n$ to $n-2m-1$ in the second sum of \eqref{h1},  and then replace $m$ by $m-1$ ,  change the order of  summation. This summation becomes
\begin{equation}\label{sum-h2}
  \sum_{m=0}^{\infty}
\sum_{n=1}^{\infty}(-1)^{m+n}z^{-n}q^{{n+1\choose 2}+(2m+1)n+m^2+m}
= \sum_{n=0}^{\infty}\sum_{m=1}^{[n/2]}
    (-1)^{m+n}z^{2m-n-1}q^{\frac{1}{2}(n^2-2m^2)+\frac{1}{2}n}.
\end{equation}
Plugging  \eqref{sum-h1} and \eqref{sum-h2} into \eqref{h1}, we get the desired identity \eqref{Garvan-n2}. Thus we complete the proof of Theorem \ref{Garvan-n2-th}.\qed

Now it's an easy matter to prove \eqref{Garvan2}. 

\noindent{\it Proof of Garvan's identity \eqref{Garvan2}.} By the  relation \eqref{r-h-g2}, we see that
 \begin{align}\label{r-h-g2-pf}
 &(1+z)(zq, z^{-1}q,q;q)_{\infty}H(z;q)\nonumber \\[3pt]
 &=(1-z)(zq, z^{-1}q,q;q)_{\infty}+2z(1-z)(zq, z^{-1}q,q;q)_{\infty}g_2(z;q).
 \end{align}
Invoking    Jacobi's  triple product identity \eqref{Jacobi-o}
 and  Theorem \ref{Garvan-n2-th} in \eqref{r-h-g2-pf}, we find that
\begin{align*}
&(1+z)(zq, z^{-1}q,q;q)_{\infty}H(z;q)\\
&=\sum_{n=0}^{\infty}(-1)^nz^{-n}q^{n+1\choose2}+\sum_{n=1}^\infty
(-1)^nz^{n}q^{n\choose2}+2\sum_{n=0}^{\infty}\sum_{m=0}^{[n/2]}
    (-1)^{m+n}z^{n-2m+1}q^{\frac{1}{2}(n^2-2m^2)+\frac{1}{2}n}\\
   &\quad+2\sum_{n=0}^{\infty}\sum_{m=1}^{[n/2]}
    (-1)^{m+n}z^{2m-n}q^{\frac{1}{2}(n^2-2m^2)+\frac{1}{2}n}.
\end{align*}
Replacing $n$ by $n+1$ in the second term, and subtracting the terms with $m=0$ in the
third sum, we  obtain
\begin{align*}
  &(1+z)(zq, z^{-1}q,q;q)_{\infty}H(z;q)\\
  &= \sum_{n=0}^{\infty}(-1)^nz^{-n}q^{n+1\choose2}-\sum_{n=0}^\infty
(-1)^nz^{n+1}q^{n+1\choose2}+2\sum_{n=0}^{\infty}(-1)^nz^{n+1}q^{n+1\choose2}\\ &\quad+2\sum_{n=0}^{\infty}\sum_{m=1}^{[n/2]}
    (-1)^{m+n}(z^{n-2m+1}+z^{2m-n})q^{\frac{1}{2}(n^2-2m^2)+\frac{1}{2}n}\\
  &=\sum_{n=0}^{\infty}(-1)^n(z^{n+1}+z^{-n})q^{n+1\choose2} +2\sum_{n=0}^{\infty}\sum_{m=1}^{[n/2]}
    (-1)^{m+n}(z^{n-2m+1}+z^{2m-n})q^{\frac{1}{2}(n^2-2m^2)+\frac{1}{2}n}\\
    &=\sum_{n=0}^{\infty}\sum_{\mid m\mid\leqslant[n/2]}
  (-1)^{n+m}(z^{n-2|m|+1}+z^{2|m|-n})
  q^{\frac{1}{2}(n^2-2m^2)+\frac{1}{2}n}.
\end{align*}
Thus we obtain Garvan's identity  \eqref{Garvan2}.\qed

We proceed to prove Theorem \ref{gener2-thm}.

\noindent{\it Proof of Theorem \ref{gener2-thm}.} Substituting the Bailey pair
 in lemma \ref{Garvan-BL-2} with the conjugate Bailey pair in Lemma \ref{lem1} with $a=z^{-1}q$, and $b=zq$ into the Bailey transform, we obtain
\begin{align}
  \sum_{n=0}^{\infty}\beta_n^{(k)}\delta_n&=(zq, z^{-1}q,q;q)_{\infty}\sum_{n=0}^{\infty}
\frac{(q^2;q)_{2n}q^n}{(zq;q)_n(z^{-1}q;q)_{n+1}}\nonumber\\
&\quad\ \ \ \times\sum_{n_1,\cdots,n_{k-1}=0}^{\infty}
\frac{q^{N_1^2+\cdots+N_{k-2}^2+\frac{3}{2}N_{k-1}^2+2(N_1+\cdots+N_{k-2})
+\frac{3}{2}N_{k-1}}(-q;q)_{n_{k-1}}}
{(q;q)_{n-N_1}(q;q)_{n_1}\cdots(q;q)_{n_{k-2}}(q^2;q)_{2n_{k-1}}}\nonumber\\
&=\sum_{m=0}^{\infty}\alpha_{2m}^{(k)}\gamma_{2m}\nonumber\\
&=\sum_{m=0}^{\infty}(-1)^mq^{(4k-3)(m^2+m)}\left(1+\sum_{n=1}^{\infty}
(-1)^nq^{n+1\choose 2}((z^{-1}q^{2m+1})^n+(zq^{2m})^n)\right)\nonumber 
\end{align}
\begin{align}
&\hskip -2cm =\sum_{m=0}^{\infty}\sum_{n=0}^{\infty}
(-1)^{m+n}z^{n}q^{{n+1\choose 2}+2mn+(4k-3)(m^2+m)}\nonumber\\\label{h2}
&\hskip -2cm \quad\ \ +\sum_{m=0}^{\infty}
\sum_{n=1}^{\infty}(-1)^{m+n}z^{-n}q^{{n+1\choose 2}+(2m+1)n+(4k-3)(m^2+m)}.
\end{align}
Replacing $n$ by $n-2m$ in the first term of \eqref{h2}  and interchanging the order of summation,
 we find that
\begin{align}\label{sum-h3}
& \sum_{m=0}^{\infty}\sum_{n=0}^{\infty}
(-1)^{m+n}z^{n}q^{{n+1\choose 2}+2mn+(4k-3)(m^2+m)} \nonumber\\
&\quad = \sum_{n=0}^{\infty}\sum_{m=0}^{[n/2]}
(-1)^{m+n}z^{n-2m}q^{\frac{1}{2}(n^2-2m^2)+\frac{1}{2}n+4(k-1)(m^2+m)}.
\end{align}
Replace $n$ by $n-2m-1$ in the second sum of \eqref{h2}, and then shift $m$ to $m-1$,  reverse the order of summation to yield
\begin{align}\label{sum-h4}
 &  \sum_{m=0}^{\infty}\sum_{n=1}^{\infty}(-1)^{m+n}z^{-n}q^{{n+1\choose 2}+(2m+1)n+(4k-3)(m^2+m)}\nonumber \\
 &\quad  =\sum_{n=0}^\infty\sum_{m=1}^{[n/2]}(-1)^{m+n}z^{2m-n-1}
q^{\frac{1}{2}(n^2-2m^2)+\frac{1}{2}n+4(k-1)(m^2-m)}.
\end{align}
Substituting \eqref{sum-h3} and \eqref{sum-h4} into \eqref{h2}, we get the desired identity \eqref{gener2}. Thus, we show Theorem \ref{gener2-thm}.\qed

\section{A two-variable Hecke-Rogers identity for $K(z;q)$ and generalization}

In this section, we give a proof of Garvan's identity \eqref{Garvan4} by the method of Bailey pair. We then give a generalization  \eqref{gener3} of this identity.

\medskip

\noindent{\it Proof of Garvan's identity \eqref{Garvan4}.}
Applying the Bailey transform with the conjugate Bailey pair in Lemma
\ref{lem1} by setting $a=z$ and $b=z^{-1}q$ and the Bailey pair in Lemma \ref{Garvan-BL-3}, we find that
\begin{align}
   \sum_{n=0}^{\infty}{\beta_n\delta_n}
   &=(zq,z^{-1}q,q;q)_{\infty}
   \sum_{n=0}^{\infty}\frac{(-1)^nq^{n^2/2}(q^{1/2};q)_n}
   {(zq;q)_n(z^{-1}q;q)_n}\nonumber\\
  &= \sum_{n=0}^{\infty}{\alpha_n\gamma_n}\nonumber\\
    &=\sum_{m=0}^{\infty}(-1)^mq^{(m^2+m)/4}\left(1+\sum_{n=1}^{\infty}
    (-1)^nq^{n\choose 2}((zq^{m+1})^n+(z^{-1}q^{m+1})^n)\right)\nonumber
     \end{align}
\begin{align}\label{k1}
    &\hskip -4cm =\sum_{m=0}^{\infty}\sum_{n=0}^{\infty}(-1)^{m+n}z^{-n}
    q^{{n+1\choose 2}+\frac{m^2+m}{4}+mn}\nonumber \\[2pt]
    &\hskip -4cm \quad +\sum_{m=0}^{\infty}\sum_{n=1}^{\infty}
    (-1)^{m+n}z^nq^{{n+1\choose 2}+\frac{m^2+m}{4}+mn}.
    \end{align}
Shifting $n$ to $n-m$ in the first sum of \eqref{k1}, and interchanging the order of summation,  we  get
    \begin{equation}\label{sum-k1}
      \sum_{m=0}^{\infty}\sum_{n=0}^{\infty}(-1)^{m+n}z^{-n}
    q^{{n+1\choose 2}+\frac{m^2+m}{4}+mn}
    = \sum_{n=0}^{\infty}\sum_{m=0}^n(-1)^nz^{m-n}
  q^{\frac{1}{4}(2n^2-m^2)+\frac{1}{4}(2n-m)}.
    \end{equation}
   For the second sum of \eqref{k1},  we first shift  $n$ to $n-m$,   then    replace $m$ by $m-1$, and  reverse the order of summation, thus  the second sum of \eqref{k1} becomes
    \begin{equation}\label{sum-k2}
        \sum_{m=0}^{\infty}\sum_{n=1}^{\infty}
    (-1)^{m+n}z^nq^{{n+1\choose 2}+\frac{m^2+m}{4}+mn}
    =\sum_{n=0}^{\infty}\sum_{m=1}^n(-1)^nz^{n-m+1}
    q^{\frac{1}{4}(2n^2-m^2)
  +\frac{1}{4}(2n+m)}.
    \end{equation}
Substituting  \eqref{sum-k1} and \eqref{sum-k2} into \eqref{k1} and replacing $q$ by $q^2$,  we  recover     \eqref{Garvan4}. This completes the proof. \qed

Using the similar argument, we could show \eqref{gener3} is valid.
\medskip

\noindent{\it Proof of Theorem \ref{gener3-thm}.}  Applying the Bailey transform with Bailey pair in Lemma \ref{Garvan-BL-3-gen} and the conjugate Bailey pair in Lemma \ref{lem1} with $a=z$, and $b=z^{-1}q$,   we have
\begin{align}
\sum_{n=0}^{\infty}\beta_n^{(k)}\delta_n
&=(zq, z^{-1}q,q;q)_{\infty}\sum_{n=0}^{\infty}
\frac{(q;q)_{2n}q^n}{(zq;q)_n(z^{-1}q;q)_{n}}\nonumber\\
&\quad \times\sum_{n_1,\cdots,n_{k-1}=0}^{\infty}
\frac{(-1)^{n_{k-1}}q^{N_1^2+\cdots+N_{k-2}^2+\frac{3}{2}N_{k-1}^2
+N_1+\cdots+N_{k-2}}(q^{\frac{1}{2}};q)_{n_{k-1}}}
{(q;q)_{n-N_1}(q;q)_{n_1}\cdots(q;q)_{n_{k-2}}
(q;q)_{2n_{k-1}}}\nonumber\\
&=\sum_{n=0}^{\infty}\alpha_n^{(k)}\gamma_n\nonumber\\
&=\sum_{m=0}^{\infty}(-1)^mq^{\frac{4k-3}{4}(m^2+m)}
\left(1+\sum_{n=1}^{\infty}(-1)^nq^{n+1\choose 2}
((zq^m)^n+(z^{-1}q^m)^n)\right)\nonumber\\
&=\sum_{m=0}^{\infty}\sum_{n=0}^{\infty}(-1)^{m+n}z^{-n}
q^{{n+1\choose 2}+mn+\frac{4k-3}{4}(m^2+m)}\nonumber\\\label{k2}
&\quad+\sum_{m=0}^{\infty}\sum_{n=1}^{\infty}(-1)^{m+n}z^{n}
q^{{n+1\choose 2}+mn+\frac{4k-3}{4}(m^2+m)}.
\end{align}
Replace $n$ by $n-m$ in the first sum of \eqref{k2} and interchange the order of summation to get
\begin{align}\label{sum-k3}
 &\sum_{m=0}^{\infty}\sum_{n=0}^{\infty}(-1)^{m+n}z^{-n}
q^{{n+1\choose 2}+mn+\frac{4k-3}{4}(m^2+m)}\nonumber \\[2pt]
&\quad = \sum_{n=0}^{\infty}\sum_{m=0}^n(-1)^{n}z^{m-n}
q^{\frac{1}{4}(2n^2-m^2)+\frac{1}{4}(2n-m)+(k-1)(m^2+m)}.
\end{align}
Shifting $n$ to $n-m$ in the second sum of \eqref{k2},   replacing  $m$ by $m-1$,  and then reversing the order of summation, this gives
\begin{align}\label{sum-k4}
&\sum_{m=0}^{\infty}\sum_{n=1}^{\infty}(-1)^{m+n}z^{n}
q^{{n+1\choose 2}+mn+\frac{4k-3}{4}(m^2+m)}\nonumber \\[2pt]
&\quad =\sum_{n=0}^{\infty}\sum_{m=1}^n(-1)^{n}z^{n-m+1}
q^{\frac{1}{4}(2n^2-m^2)+\frac{1}{4}(2n+m)+(k-1)(m^2-m)}.
\end{align}
Hence  \eqref{gener3} is obtained  upon substituting \eqref{sum-k3} and \eqref{sum-k4} in \eqref{k2} and replacing $q$ by $q^2$. This completes the proof. \qed

\noindent{\bf Acknowledgments.} This work was supported by the 973 Project, the PCSIRT Project of the Ministry of Education and the National Science Foundation of China.

\end{document}